\documentclass{article}
\usepackage{amsmath,amsthm,amsfonts}

\input epsf.tex
\newdimen\epsfxsize
\newdimen\epsfysize
\def\qed{\vrule height5pt width3pt depth.5pt}

\theoremstyle{plain}
\newtheorem{thm}{Theorem}[section]
\newtheorem{cor}[thm]{Corollary}
\newtheorem{lem}[thm]{Lemma}
\newtheorem{prop}[thm]{Proposition}
\newtheorem{conj}[thm]{Conjecture}

\newtheorem{rem}{Remark}[section]

\begin{document}

\markboth{H. A. Dye and Louis H. Kauffman}
{Virtual Homotopy}

\title{Virtual Homotopy}

\author{H. A. Dye\\ McKendree University \\
hadye@mckendree.edu\\ Louis H. Kauffman \\ Department of Mathematics, Statistics, \& Computer Science \\
  University of Illinois at Chicago\\
 kauffman@uic.edu} 

\maketitle


\begin{abstract} Two welded (respectively virtual) link diagrams are homotopic if one may be transformed into the other by a sequence of extended Reidemeister moves, classical Reidemeister moves, and self crossing changes. In this paper, we extend Milnor's $ \mu $ and $ \bar{ \mu } $ invariants to welded and virtual links. We conclude this paper with several examples, and compute the $ \mu $ invariants using the Magnus expansion and Polyak's skein relation for the $ \mu $ invariants.
\end{abstract}

\maketitle

\section{Introduction}

\begin{figure}[htb] \epsfysize = 1 in
\centerline{\epsffile{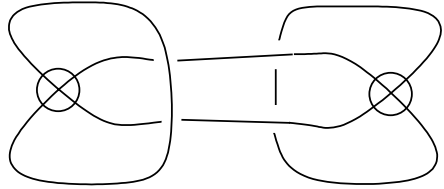}}
\caption{Kishino's knot}
\label{fig:kish}
\end{figure}

Virtual link diagrams are a generalization of classical link diagrams that was introduced by Louis H. Kauffman in 1996 \cite{kvirt}. A virtual link diagram is a decorated immersion of $ n $ copies of $ S^1 $ with two types of crossings: classical and virtual. Classical crossings are indicated by over/under markings and virtual crossings are indicated by a solid encircled X. An example of a virtual link diagram is shown in figure \ref{fig:kish}.

\begin{figure}[htb] \epsfysize = 0.5 in
\centerline{\epsffile{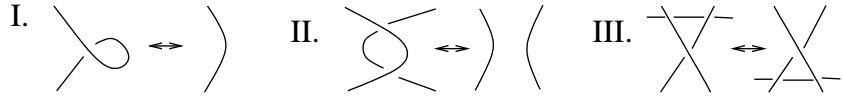}}
\caption{Classical Reidemeister moves}
\label{fig:rmoves}
\end{figure}

\begin{figure}[htb] \epsfysize = 1.25 in
\centerline{\epsffile{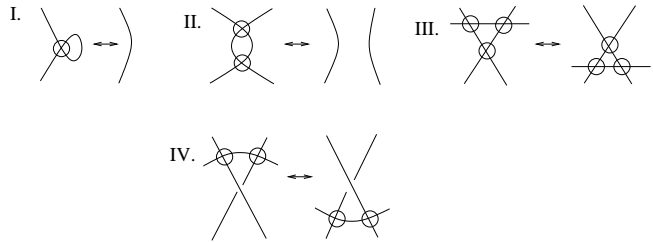}}
\caption{The virtual Reidemeister moves}
\label{fig:vrmoves}
\end{figure}

\begin{figure}[htb] \epsfysize = 0.75 in
\centerline{\epsffile{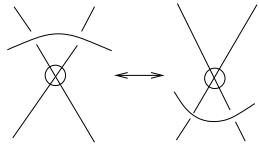}}
\caption{The upper forbidden move}
\label{fig:forbidmove}
\end{figure}

Two virtual link diagrams are said to be \emph{virtually isotopic} if one may be transformed into another by a sequence of classical Reidemeister moves (shown in \ref{fig:rmoves}) and the virtual Reidemeister moves shown in figure \ref{fig:vrmoves}. Two virtual link diagrams are said to be \emph{welded equivalent} if one may be transformed in the other by a sequence of Reidemeister moves (classical and virtual) and the upper forbidden move (shown in figure \ref{fig:forbidmove}). When we use welded equivalence, we shall refer to the virtual diagrams as welded diagrams. The forbidden move is not allowed during virtual isotopy.

In this paper, we apply a variety of results from the area of virtual knot theory. 
We will extend the definition of link homotopy to virtual and welded link diagrams, and Milnor's $ \bar{\mu} $ invariants for virtual and welded links. 

We define two virtual link diagrams to be \emph{virtually homotopic} if one diagram may be transformed into the other by virtual isotopy and self crossing change. (By self crossing change, we mean changing the over/under markings at a crossing between two segments of the same link component.) Two welded diagrams are \emph{welded homotopic} if one diagram may be transformed into the other by a sequence of classical and virtual Reidemeister moves, the forbidden move, and self crossing change.

Homotopic classical link diagrams (and braids) have been studied extensively, starting with Milnor \cite{milnor1} \cite{milnor2}, Goldsmith \cite{goldsmith}, and more recently Habegger and Lin \cite{xsl}. There are significant differences between the study of homotopy classes in the classical case and the virtual case. We introduce these differences by examining the homotopy class of several virtual and welded knot diagrams.

\begin{figure}[htb] \epsfysize = 0.55 in
\centerline{\epsffile{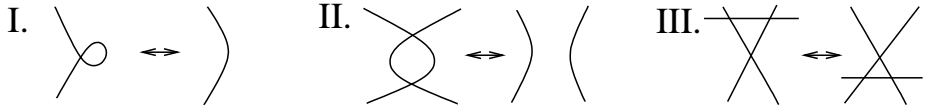}}
\caption{Flat Reidemeister moves}
\label{fig:flatmoves}
\end{figure}

\begin{figure}[htb] \epsfysize = 0.7 in
\centerline{\epsffile{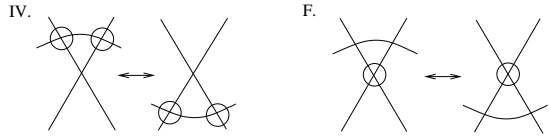}}
\caption{The flat virtual Reidemeister move and the flat forbidden move}
\label{fig:flatvrmoves}
\end{figure}

We first consider virtual knot and link diagrams.  
To distinguish homotopy classes of virtual knot diagrams, we recall the category of  \emph{flat virtual knot diagrams}. A flat virtual knot diagram is a decorated immersion of $ S^1 $ into the plane with two types of crossings: flat crossings and virtual crossings. (Flat crossings are indicated by an unmarked $ X $.) Two flat diagrams are equivalent if one may be transformed into another by a sequence of flat Reidemeister moves and virtual Reidemeister moves. The flat classical Reidemeister moves are shown in figure \ref{fig:flatmoves}. The flat versions of virtual move IV. and the forbidden move are shown in figure \ref{fig:flatvrmoves}.

\begin{figure}[htb] \epsfysize = 0.5 in
\centerline{\epsffile{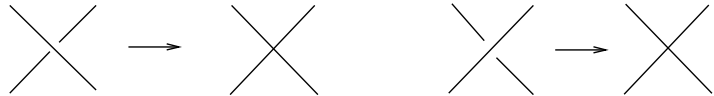}}
\caption{Mapping into the flat knot diagrams}
\label{fig:flatten}
\end{figure}

We map the category of virtual knot diagrams into the category of flat virtual knot diagrams by \emph{flattening} all the classical crossings in the diagram as indicated in figure 
\ref{fig:flatten}.

We observe that equivalence classes of homotopic virtual knot diagrams are in 1-1 correspondence with equivalence classes of flat virtual knot diagrams. The existence of non-trivial flat virtual diagrams demonstrates that not all virtual knot diagrams are homotopy equivalent to the unknot. (Recall that all classical knot diagrams are homotopic to the unknot.) 
In fact, we prove the following theorem:
\begin{thm} The category of flat virtual link diagrams is an invariant image of the category of homotopic virtual link diagrams. In the case of knot diagrams, the category of flat virtual knot diagrams is equivalent to the category of homotopic virtual link diagrams. 
\end{thm}
\textbf{Proof:} 
Let $ L $ and $L'$ be homotopic virtual link diagrams. Then there is a finite sequence of Reidemeister moves and self crossing changes transforming $L $ into the $L'$. The $ Flat $ version  of this sequence may be applied to $ Flat(L) $. This demonstrates that the category of flat virtual links is an invariant image of the category of homotopic virtual link diagrams.

We restrict our attention to knot diagrams. Let $ K $ and $K'$ be virtual knot diagrams. Suppose that $ Flat(K) $ is equivalent to $Flat(K')$. Then there exists a sequence of flat and virtual Reidemeister moves transforming $ Flat(K)$ into $Flat(K')$. Each step of the sequence may be realized as a virtual or classical Reidemeister move by switching the classical crossings appropriately. Hence $ K $ is homotopic to the $K'$.\qed 

We immediately obtain the following corollary.
\begin{cor} A virtual knot diagram $K $ is homotopic to the unknot if and only if $ Flat(K) $ is equivalent to the unknot.
\end{cor}
Note that as a result, Kishino's knot is not equivalent to the unknot.  It was shown in \cite{kadokami} that there is no sequence of flat moves that transforms a flat Kishino's knot into the unknot. 

\begin{thm} If $ K $ is a welded knot then $K $ is homotopic to the unknot. 
\end{thm}
\textbf{Proof:}Two welded knots are equivalent if one can be transformed into the other by a sequence of virtual and classical moves, including the forbidden move. Note that in the forbidden move, two classical crossings pass over a virtual crossing. We can switch the crossings on this move, so that two classical crossing pass under a virtual crossing. Using both the under and over version of the forbidden move, any virtual knot can be unknotted \cite{fop}.\qed 
 
\smallskip

The category of flat virtual links (more than one component) is not equivalent to homotopy classes of virtual links. We can not flatten crossings between the components. Similarly, the approach using both forbidden moves can not be applied to welded link diagrams. To unknot and unlink a link diagram, we may need  to switch crossings between two different components of the link diagram. In the next section, we generalize  Milnor's link group to virtual and welded link diagrams.

\section{Link Groups}

\begin{figure}[htb] \epsfysize = 1 in
\centerline{\epsffile{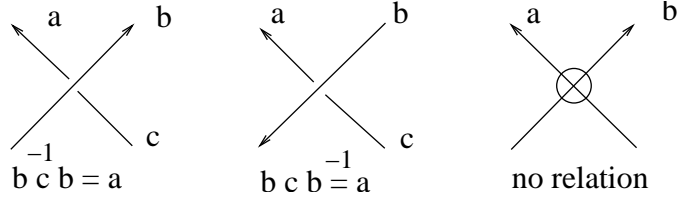}}
\caption{Fundamental group relations}
\label{fig:fund}
\end{figure}

In \cite{milnor1}, Milnor defined the \emph{link group} of a classical link. We extend this definition to virtual and welded link diagrams. 
\begin{rem} The link group is also referred to as the reduced group in \cite{xsl} and \cite{milnor1}. 
\end{rem}
Let $ L $ be a  virtual link diagram with $ n $ components: $ K_1, K_2 , \ldots K_n$. Recall that the fundamental group of the virtual link digram, $ \pi_1 (L) $, is a free group modulo relations determined by the classical crossings.
Generators of the fundamental group correspond to arcs (that initiate and terminate at a classical crossing). The arcs of a particular component are referred to as \emph{meridians} belonging to that component. 
The labels: $ a_{i1}, a_{i2}, \ldots a_{ik} $, are the arcs of the $ i^{th} $ component. We then record relations, $ r_1, r_2, \ldots r_m $, that are determined by each crossing, as shown in figure \ref{fig:fund}.

The fundamental group of a link $L $ is:
\begin{equation*}
 \pi_1 (L) = ( a_{11}, a_{12}, \ldots a_{1k}, a_{21}, \dots a_{n\hat{k}} | r_1, r_2, \ldots r_m )
\end{equation*}

Assume that the crossings in the diagram have been labeled with the integers $ 1,2, \ldots N $. The longitude of a component is defined to be the product of the generators encountered at successive underpasses in the opposite direction of the component's orientation. We can then write the longitude of a component as $ l= l_{i_1} ^{ \epsilon_1} l_{i_2} ^{ \epsilon_2} \ldots l_{i_n} ^{\epsilon_n} $ where $l_{i_k} $ (equal to some $a_{kr}$) is the label of the the overcrossing at the $i_k ^{th} $ crossing underpassed. If the crossing is positively oriented then $ \epsilon_k =1 $, otherwise $\epsilon_k = -1 $. Note that the longitude is well defined up to cyclic permutation of its factors. 

\begin{rem}
This definition of longitude differs from the definition of longitude given in \cite{kvirt}. In \cite{kvirt}, the longitude is defined to be the product of the generators encountered in the direction of the orientation. Up to the direction of traversal, the longitude agrees
with the usual geometric notion of longitude for fundamental group of a
classical knot.   We use this definition to construct the $ \mu $ and $ \bar{ \mu } $ invariants based on Artin's original definition \cite{milnor1}. 
\end{rem}

We recall that for two elements $a$ and $b$ of a group $G$, the commutator has the form: $ a b a^{-1} b^{-1} $. We denote 
the commutator of $ a $ and $b$ as:
\begin{equation} \label{commutat}
[a,b].
\end{equation}
The commutators of the generators for the $ i^{th} $ component have the form: $ a_{ij} a_{ik} a_{ij} ^{-1} a_{ik} ^{-1} $ and are denoted $ [ a_{ij},a_{ik}] $. Let $ A $ to denote the smallest normal subgroup generated by the commutators formed by the meridians of a fixed component: $ [a_{ij},b_{ik}] $ from all components ($ i \in \lbrace 1, 2 , \dots n \rbrace $) of the diagram.
The \textit{link group} of $ L $, denoted as $ G(L) $, is the fundamental group, $ \pi_1 (L) $, modulo $ A $:
\begin{equation*}
G(L)= \pi_1 (L) / A ).
 \end{equation*}
Milnor proved the following theorem.
\begin{thm}Let $ L $ be a classical link diagram and let $ \hat{L} $ be the trivial link. If $ G(L)  $ is isomorphic to $ G(\hat{L} ) $, via an isomorphism preserving the conjugacy classes corresponding to meridians, then $ L $ is homotopic to the trivial link .
\end{thm}
\textbf{Proof:} See \cite{milnor1}.

For virtual link diagrams, the requirement that $ G(L) \cong G(\hat{L}) $
 is necessary but not sufficient.  However, we can prove the following theorem using the Wirtinger presentation and extend the result to virtual link diagrams.

\begin{figure}[htb] \epsfysize = 1 in
\centerline{\epsffile{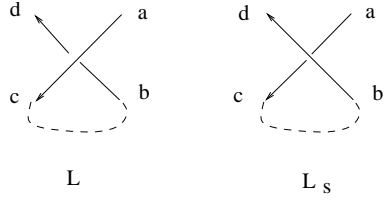}}
\caption{Invariance under crossing change}
\label{fig:crosschange}
\end{figure}

\begin{thm} If two virtual (welded) link diagrams, $ L $ and $ \hat{L} $, are homotopic then $ G(L)$ is isomorphic to $G( \hat{L} ) $.
\end{thm}
\textbf{Proof: } The fundamental group is invariant under the Reidemeister moves. We need only prove that the link group is invariant under crossing change. We examine the case when two diagrams differ by exactly one crossing, as illustrated in figure \ref{fig:crosschange}. 
In this figure, the dotted line indicates that both segments of the crossing are in the same component. 
(Alternate choices for this connection are possible, but the arguments are similar to the one 
given below.)
From the crossing, $L $, in figure \ref{fig:crosschange} we obtain the following relationships:
\begin{equation} \label{little} 
a=c
\end{equation}
\begin{equation} \label{little2}
bc=ad.
\end{equation}
Since the segments $ b $ and $ c $ are part of the same component, we observe that
\begin{equation*}
 b = g c g^{-1} \text{ where }  g \in \pi_1 (L) .
 \end{equation*}
Substituting equation \ref{little} into $ bc =ad $, 
we determine that $ b=c d c^{-1} $. But $ d $ and $ c $ are meridians of the same component, so that $[c,d] \equiv 1 $ in $G(L)$. As a result, we observe that $ b \equiv d $ in $ G(L) $.
In the diagram, $ L_s $, from figure \ref{fig:crosschange}:
\begin{gather*}
b=d \\
d a = b c
\end{gather*}
Rewriting the last equation, we observe that
$ b c b^{-1} =a $. Note that $ b $ and $c $ are meridians of the same component, so that $ c \equiv a $ in $ G(L_s) $. The two diagrams only differ at a single crossing and we 
 conclude that $ G(L) $ and $ G(L_s) $ are isomorphic.\qed

Many non-trivial virtual link diagrams have a fundamental group equivalent to the fundamental group of the unknot; the link group does not determine the virtual homotopy class.
For example, Kishino's knot has a trivial link group (since $ \pi_1 (K) = \mathbb{Z} $) but Kishino's knot is not virtually homotopic to the unknot. However,  Kishino's knot is welded homotopic to the unknot.

For welded link diagrams we make the following conjectures:

\begin{conj} If $L$ is welded link diagram and $ \pi_1 (L) \cong \mathbb{Z} \oplus \ldots \oplus \mathbb{Z} $ then $ L $ is welded equivalent to the unlink.  
\end{conj}

\begin{conj} If $ L $ is a welded link diagram  and $G(L) \cong \mathbb{Z} \oplus \ldots \oplus \mathbb{Z} $ then $ L $ is homotopic to the unlink.
\end{conj}

We now extend Milnor's $ \bar{ \mu } $ invariants by utilizing Milnor's original definition of the invariant.

\section{The $ \mu $ Invariants}

We extend Milnor's $ \mu $ and $ \bar{ \mu } $ invariants to $k $ component virtual (and welded) link diagrams by
computing the representatives of the longitude of each component in the Chen group: $ \pi_1 (L) / \pi_1  (L)_{n} $. The Chen group and the group $ \pi_{1}(L)_{n}$ will be defined below. We then apply the Magnus expansion to to the longitude of each component. For each component $ i $, we obtain a polynomial, $P_i $, in $k $ variables. We then compute the modulus of each polynomial, producing the $ \bar{ \mu } $ invariants. (The classical analog of this proof is given in  \cite{milnor2}.)   Finally, we will show that the $ \bar{ \mu } $ invariants are unchanged by virtual link homotopy.

Let $ \hat{L} $ be an $k$ component virtual (or welded) braid diagram with all strands oriented downwards such that the closure of $ \hat{L} $ is the link diagram $L $.  We show that the $n^{th} $ Chen group of a $k$-component link is isomorphic to a group with $k$ generators following \cite{fenn}. 

Arrange the link $L$ so that it is the closure of a $k$ component tangle. Label the uppermost meridian of each strand as 
$ a_{10}, a_{20}, a_{30}, \ldots a_{k0} $.  After each undercrossing, continue labeling so that meridional component after $a_{ij}$ is denoted $ a_{ij+1}$. The second index in $a_{ij} $ counts the number of underpasses from $a_{i0} $ to $ a_{ij} $. The terminal meridian of each strand is labeled as $ a_{1t_1}, a_{2t_2}, a_{3t_3}, \ldots a_{kt_k} $.

\begin{figure}[htb] \epsfysize = 2.0 in
\centerline{\epsffile{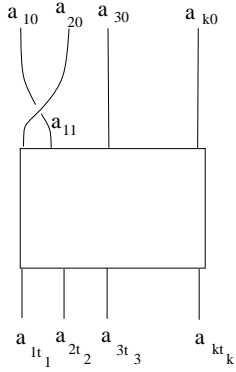}}
\caption{Labeled $k$ component tangle}
\label{fig:ktangle}
\end{figure}
Using the 
Wirtinger presentation (figure \ref{fig:fund}), we compute the relation at each classical crossing. Let $r_{ij} $ denote the overcrossing strand between $a_{ij} $ and $a_{ij+1}$ (the $j^{th}$ and $j+1 ^{th} $ segments on the $i^{th} $ component of the braid. We compute the first relation determined by the $i^{th}$ strand:
\begin{equation*}
r_{i1} ^{\epsilon_1} a_{i0} r_{i1} ^{-\epsilon_1} = a_{i1} .
\end{equation*}
In general, we find that the $j^{th}$ overcrossing of the $i^{th}$ strand determines the following relation:
\begin{equation*}
 r_{ij} ^{\epsilon_j}  a_{i(j-1)} r_{ij}^{-\epsilon_j} = a_{ij} .
\end{equation*}
Let $l_{ij} $ denote $ r_{ij} ^{ \epsilon_j} \ldots r_{i2} ^{\epsilon_2} r_{i1} ^{\epsilon_1} $.
Then $l_{ij} $ is word consisting of the $ a_{ij} $'s and
\begin{equation*}
a_{ij} = l_{ij} a_{i0} l_{ij} ^{-1}.
\end{equation*}
Since we take the closure of  $k-k$ tangle shown in figure \ref{fig:ktangle} to form the link diagram, we observe that 
$a_{it_i} = a_{i0} $ and therefore
\begin{equation*}
l_{it_i} a_{i0} l_{it_i} ^{-1} = a_{i0}. 
\end{equation*}
Note that $l_{it_i} $ is the longitude of the $i^{th}$ component.
Recall the definition of a commutator (equation \ref{commutat}) and we obtain the following presentation of the fundamental group:
\begin{equation*}
\pi_1 (L) = ( a_{ij} | [l_{it_i}, a_{i0}]=1 \text{ and }  l_{ij} a_{i0} l_{ij} ^{-1} = a_{ij} \text{, }  1 \leq i \leq k  \text{ and } 0 \leq j < t_i  ).
\end{equation*} 
\begin{rem} In this proof, we do not directly consider virtual crossings. These crossings do not contribute to the fundamental group. However, the definition of the ideal $D_{i,n} $ (which appears later in the proof) is subtly different from the classical case. \end{rem}

We define the \textit{ $n^{th}$ commutator subgroup} of a group $A$. Let $ A_1 = A $ then $ A_n = [A,A_{n-1} ]$. We use this definition in a lemma and in the definition of the Chen group.

\begin{lem} \label{com} Let $V$ and $W$ be elements of the group $A$. Let $ V \equiv \hat{W} $ modulo $A_{n-1} $ then
\begin{equation*}
V b V^{-1} \equiv W b W^{-1} \text{ modulo } A_n. 
\end{equation*}
\end{lem}
\textbf{Proof:}
Let $ V \equiv W $ modulo $ A_{n-1} $. Then $ W^{-1} V $ is an element of $ A_{n-1}$. Therefore,
$ W^{-1} V b V^{-1} W b^{-1} $ is an element of $A_n$. Hence $ W b W^{-1} \equiv  V b V^{-1} $ modulo $ A_n$.\qed

Let $ F $ denote $ \pi_1 (L) $ and let $F_n$ denote the $n^{th} $ commutator subgroup. \textit{ We denote the $n^{th} $ Chen group of $ \pi_1 (L)$  as $ \pi_1 (L)_n $ where $ \pi_1 (L)_n = F/F_n$.} 
Let $ P $ denote the free group on $k$ elements generated by $ \lbrace B_1, B_2, \ldots B_k \rbrace $. Similiarly, let $ P_n $ denote the $n^{th} $ commutator subgroup of $P$. By placing a set in the list of relations for a presented group, we mean that each element of this set becomes a relation $x=1$ for an element $x$ of that set.
We define the group $ G(n) $ as a quotient of $P$:
\begin{equation}
G(n) = \lbrace B_1, B_2, \ldots B_k | P_n, [ \rho^n (l_{it_i} ), B_i ] = 1 \text{ for } 1 \leq i \leq k \rbrace.
\end{equation}
where $ \rho^n(l_{it_i} ) $ is the image of the longitude of the $i^{th} $ component from the homomorphism
$ \rho^n : \pi_1 (L)_n \rightarrow G(n) $, which we define below and show to be an isomorphism.
 
We define a sequence of  homomorphisms $ \rho^n : \pi_1 (L)_n \rightarrow G(n) $ inductively.  
\begin{gather*}
\rho^1 (a_{ij}) = 1  \\
\rho^2 (a_{ij}) = B_i  \\
\rho^n (a_{i0}) = B_i \text{ for all } n \geq 2 \\
\rho^n (a_{ij}) \equiv  \rho^{n-1} (l_{ij} ) B_i \rho^{n-1} (l_{ij} ^{-1} ) \text{ modulo } P_{n-1}
\end{gather*}
These definitions are written in terms of the free group modulo appropriate subgroups.
Notice that  $ a_{ij} = l_{ij} a_{i0} l_{ij} ^{-1}$. 
\begin{thm}The map $ \rho^n : \pi_1 (L)_n \rightarrow G(n) $ is an isomorphism.
\end{thm}
\textbf{Proof:}
We show that $ \rho^n $ is well defined. Notice that $ G(1) $ and $ \pi_1 (L)_1 $ are both the trivial group, hence $ \rho^1 $ is well defined. For $n \geq 2$, we have made the definition
\begin{equation}
\rho^n (a_{ij}) \equiv \rho^{n-1} (l_{ij}) B_i \rho^{n-1} (l_{ij} ^{-1} ) \text{ modulo } P_{n-1}.
\end{equation}
We also observe that since $ \rho^n $ is a homomorphism then
\begin{equation}
\rho^n (a_{ij}) \equiv \rho^n (l_{ij}) B_i \rho^n (l_{ij} ^{-1}).
\end{equation}
By Lemma \ref{com}, if $ \rho^n (l_{ij} ) \equiv \rho^{n-1} ( l_{ij} ) $ modulo $ P_{n-1}$ then $ \rho^{n-1} (l_{ij}) B_i \rho^{n-1} (l_{ij} ^{-1} ) \equiv \rho^n (l_{ij}) B_i \rho^n (l_{ij} ^{-1})$  modulo $ P_{n} $. We observe that $l_{ij} $ is a product of $a_{kl} $'s so that $ \rho^n (l_{ij}) \equiv \rho^{n-1} (l_{ij}) $ modulo $ P_{n-1} $ as required.

We must also show that if $ A \equiv B $ modulo $ F_n $ then $ \rho^n (A) \equiv  \rho^n (B) $ modulo $ P_n $ for all $n \geq 2 $. 
For $ \pi_1 (L)_2 $, we consider $ [ a_{ij}, a_{kl} ] $ which is an element of $ F_2$. Note that $ \rho^2 ( a_{ij} a_{kl} a_{ij} ^{-1} ) = B_i B_k B_i $. Since $ B_i B_k B_i ^{-1} \equiv B_k $ modulo $P_2 $, this statement is true for $n=2 $.

We assume that the statement holds true for all $ \hat{n} < n $. Suppose that $ A \equiv B $ modulo $ F_k$ then there exists $W$, an element of $ F_{k-1} $ such that $ W A W^{-1} B^{-1} $ is an element of $ F_k $. Now, $ \rho^n ( W A W^{-1}) = \rho^n (W) \rho^n (A) \rho^n (W^{-1}) $. But $ \rho^n (W) \equiv \rho^{n-1} (W) $ modulo $ P_{n-1} $ and, as a result, 
$ \rho^n (W) \equiv 1 $ modulo $ P_{n-1}$. Hence $ \rho^n (A) \equiv \rho^n (B) $ modulo $P_n$.

Now, let $M $ denote the smallest normal subgroup generated by commutators of the form $ [l_{it_i}, a_{i0}]$. We observe that if $ A \equiv B $ modulo $ M $ then $ \rho^n (A) \equiv  \rho^n (B) $ modulo $ \rho^n (M) $. However, these relations are included in the definition of $ G(n)$. This completes the proof that $ \rho^n $ is injective.

We need only show that $ \rho^n $ is onto. The homomorphism $ \rho^n $ is onto since $ \rho^n (a_{i0}) = B_i $.\qed

We now switch notation to simplify the description of $ G(n)$. We use $ a_1, a_2, \ldots a_k $ to generate the group $P$.
The proof above demonstrates that the Chen group of any $k$ component link diagram has a standard presentation.
\begin{prop}
Let $ L $ be an $k $ component link diagram. Let $l_i $ denote the longitude of the $i^{th} $ component. There is a presentation of $n^{th} $ Chen group of $L$ as:
\begin{equation}
G(n) = \lbrace a_1 , a_2 \ldots a_k | P_n, [a_i, l_i ^{ (n)} ] =1 \text{ for } 1 \leq i \leq k \rbrace. 
\end{equation}
where $ l_i ^{(n)} $ is a reduced word in $ P/P_{n} $.
\end{prop}
Note that 
$
G(1) = {1}
$, the trivial group.
The group $ G(2) $ is the abelianization of the fundamental group, since any conjugate of $ a_i \equiv a_i $ modulo $P_2 $. 

We now apply the Magnus expansion to the group $ G(n) $. The Magnus expansion is a homomorphism  $ \psi $ such that:
\begin{equation*} 
\psi : G(n) \rightarrow \mathbb{Z}[[x_1, x_2, \dots x_k ]]
\end{equation*}
where $ \mathbb{Z}[[x_1, x_2, \dots x_k ]] $ is the formal power series ring on non-commuting variables $ x_1, x_2, \ldots x_k $.
 We define the Magnus expansion on $ G(n) $ by specifying the generators:
\begin{gather*} 
 \psi( a_i )= 1 + x_i \\
 \psi( a_i ^ {-1})= 1 -x_i + x_i ^2 - x_i ^ 3 \ldots 
\end{gather*}

Notice that $ \psi( a_i a_i^{-1})= \psi(1)= 1 $ but  $ \psi (a_i a_j) \neq \psi (a_j a_i)$ since the $ x_i $ do not commute. 
For convenience, we will denote $ l_i ^{(n)} $ as $ w_i $ for the remainder of this paper.
The image of the longitude of the $i^{th} $ component, $w_i $, can be written as a power series in $k$ variables with integer coefficients (where $k$ equals the number of link components).
\begin{equation} \psi(w_i) = 1 + \Sigma \mu (i_1, \dots i_s,w_i) x_{i_1} x_{i_2} ... x_{i_s} . 
\end{equation}

Let $ J $ be a proper subset of $ \lbrace i_1, i_2, \ldots i_s \rbrace $.
We define 
$ \Delta (i_1, i_2, \ldots i_s, w_i) $ as follows:
\begin{equation*}
  \Delta (i_1, i_2, \ldots i_s, w_i) =
   \text{the g.c.d. of } \lbrace \mu (J,w_i)| J \subset \lbrace i_1, i_2, \ldots i_s \rbrace \rbrace. 
\end{equation*}
Then $ \bar{ \mu } (i_1, i_2, \ldots i_s, w_i ) \equiv $
\begin{equation*}
 \begin{cases}
 \mu  ( i_1, i_2, \ldots i_s, w_i ) \text{ modulo } \Delta ( i_1, i_2, \ldots i_s, w_i ) & \text{ if } i  \notin \lbrace i_1, i_2, \ldots i_s \rbrace, \\
   0 &
 \text{ if } i \in \lbrace i_1, i_2, \ldots i_s \rbrace.
 \end{cases}
\end{equation*}

We now define an ideal $ D_{i, n} $ in $ \mathbb{Z} [[ x_1, x_2, \ldots x_k ]] $.
Let $ v= c x_{i_1} x_{i_2} \ldots x_{i_s} $. Then $ v $ is an element of $ D_{i,n} $ if  $ v $ meets any of the following criteria (note that $s$ is the length of the word in the $x_i$'s.) :
\begin{enumerate}
\item $| \lbrace i_1, i_2, \ldots i_s \rbrace | \geq n $
\item $ c \equiv 0 \text{ modulo } \Delta (i_1, i_2, \ldots i_s, w_i) $
\item $ i \in \lbrace i_1, i_2, \ldots i_s \rbrace $
\end{enumerate}
The ideal $D_{i,n} $ is generated by the monomials $v$ of this form.
We observe that if $ \bar{ \mu } (J, w_i) \neq 0 $ then  the cardinality of $J$ is less than $n$ ($ |J|  < n $).
\begin{rem} The definition of $ \Delta (i_1, i_2, \ldots i_s,w_i) $ in the virtual (welded) case differs from the classical case. In the classical case, if $ i \in \lbrace i_1, i_2, \ldots i_s \rbrace $ then $ \mu (i_1, i_2, \ldots i_s,w_i) =0 $. We will discuss the implications of this fact in sections \ref{examples} and \ref{links}.
\end{rem}
\noindent
For any $ w_i \in G(n) $, the Magnus expansion of $ w_i$ ($ \psi (w_i) $) modulo $ D_{i, n} $  
determines the values of the $ \bar{ \mu} $ as defined previously. That is, $ \bar{ \mu } (J,w_i)  \in \mathbb{Z}[[x_1, x_2, \ldots x_n ]] /  D_{i,n} $. We shall prove that these $ \bar{ \mu } $ values are homotopy invariants. In fact,
we need several facts about the ideal $ D_{i,n} $ in order to prove that
 $ \mathbb{Z} [[ x_1, x_2, \ldots x_k ]] $  modulo $  D_{i,n} $ is a homotopy invariant.
We first prove the following technical lemmas about $ D_{i,n} $.
\begin{lem} If $ v= c x_{i_1} x_{i_2} \ldots x_{i_s} $ is an element of $ D_{i,n} $, that is,
$ v $ is equivalent to $ 0 $ in $ \mathbb{Z} [[ x_1, x_2, \ldots x_k]] / D_{i,n} $ then 
$ v x_j $ and $ x_j v $ are elements of $ D_{i,n} $ for any $j \in {1,2, \ldots k}$.
\end{lem}

\textbf{Proof: } Let $ v = c x_{i_1} x_{i_2} \ldots x_{i_s} $. If $ j=i $ then $ v x_j $ and $ x_j v $ are elements of $ D_{i,n} $. If $ s \geq n-1 $ then $ x_j v $ and $ v x_j $ are elements of $ D_{i,n} $. 
If $ s < k-1 $, we observe that $ c \equiv 0 $ modulo
 $ \Delta (i_1, i_2, \ldots i_s, w_i) $. Note that $ \lbrace i_1, i_2, \ldots i_s \rbrace $ is a subset of $ \lbrace i_1, i_2, \ldots i_s, j \rbrace $. As a result,
  $ \Delta (i_1, i_2, \ldots i_s,j, w_i ) $ divides $ \Delta(i_1, i_2, \ldots i_s ) $. 
  Then $ \Delta (i_1, i_2, \ldots i_s,j, w_i ) $ divides $ c $ and $ v x_j  \in D_{i,n} $. 
  Similarly, we may argue that $ x_j v \in D_{i,n} $. \qed
  
We have demonstrated that $ D_{i,n} $ is a two sided ideal.

\begin{lem} \label{an} If $ m \in P_n $ then  $ \psi(m) - 1
\in D_{i,n} $.
\end{lem}
\textbf{Proof: }
To demonstrate that $ \psi (m) \in D_{i,n} $, we will show that each non-constant term has degree greater than or equal to $ n $. 
We consider $ m \in P_2 $. Without loss of generality, we assume that $ m=[a_i,a_j] $.
Now,
\begin{equation} \label{m1}
\psi (m) = (1+x_i) (1+x_j) (1-x_i + x_i ^2 - x_i ^3 \ldots) (1-x_j + x_j ^2 - x_j ^3 \ldots)
\end{equation}
Let 
\begin{align*} v_i &= (1-x_i + x_i ^2 - x_i ^3 \ldots) \\ 
 v_j &= (1-x_j + x_j ^2 - x_j ^3 \ldots) 
 \end{align*}
Rewriting equation \ref{m1}:
\begin{align*}
\psi(m) &= (1 + x_i + x_j + x_i x_j) v_i v_j \\
&= (1 + x_j v_i + x_i x_j v_i) v_j \\
&= (1 + x_j + (-x_j x_i + x_j x_i ^2  - x_j x_i ^3 \ldots ) + x_i x_j v_i) v_j \\
&= 1 + (- x_j x_i + x_j x_i ^2  - x_j x_i ^3 \ldots )v_j + x_i x_j v_i v_j.
\end{align*}
All non-constant terms have degrees greater than or equal to two. 
Let $ m \in P_{n-1} $, $ \psi ( m) = 1 + M $ and $ \psi( m ^{-1}) = 1 + \hat{M} $. By assumption, all terms in $ M $ and $ \hat{M} $ have degree greater than or equal to $ n-1 $. Note that $ (1+M)(1+\hat{M}) = (1 + \hat{M}) ( 1+ M) = 1 $. 
Now, without loss of generality, we may assume that all elements of $ P_n $ have the form:
\begin{equation*}
[m, a_i] \text{ where } m \in P_{n-1}.
\end{equation*}
We compute that
\begin{equation*}
\psi([m,a_i]) = (1+M) (1+x_i) (1+\hat{M}) (1-x_i + x_i ^2 - x_i ^ 3 \ldots )
\end{equation*}
Letting $ v_i $ denote $ 1- x_i + x_i ^2 - x_i ^3 \ldots  $, we multiply:
\begin{align*}
\psi([m,a_i])
 &= ( 1 + M +  x_i + x_i M + \hat{M} + M \hat{M} + x_i \hat{M}  + x_i M \hat{M}) v_i \\ 
&=(v_i + x_i v_i + M x_i v_i  + x_i \hat{M} v_i + M x_i \hat{M} v_i) \\
&= 1 + M x_i v_i + x_i \hat{M} v_i + M x_i \hat{M} v_i .
\end{align*}
Since each term in $ M $ and $ \hat{M} $ has degree greater than or equal to $ n-1 $, we observe that $ \psi([m,a_i]) \equiv 1 \text{ modulo } D_{i,n} $. \qed

\section{Invariance of $ \bar{ \mu } $ under Homotopy}
We first prove that $ \bar{ \mu } $ is a link invariant. The fundamental group is invariant under virtual isotopy and the upper forbidden move. As a result, $G(n) $ is also invariant under virtual isotopy and the upper forbidden move. We need only show that the image of $ G(n) $ is invariant under change of base point and invariant under multiplication by elements of $ P_{n} $ and $ [a_i,  w_i ] $. 

\begin{lem} \label{linkinvar} Let $ G(n) $ denote a presentation of the Chen group of a link diagram $L$. Let $w_i = l_i ^{(n)} $ (the longitude of component $i$ modulo $P_{n} $) and let $ a_j$ represent the image of the generator of $j^{th} $ component in $ G(n) $.
\begin{enumerate}
	\item $ \psi (a_j w_i a_j ^{-1}) = \psi (w_i) $ and $\psi(a_j a_i a_j ^{-1} ) = \psi (a_i)$ (invariance under conjugation by $ a_j $,
	 representing a change of base point).
	\item $ \psi( [a_j, w_j]) = 1$ (verification of  $ [a_j, w_j ] =1 $).
	\item $ \psi(m) = 1$ for  $ m \in P_{n}  $ (verification of the identity $ m \equiv 1 $ for $ m \in  P_n $).
\end{enumerate}	
\end{lem}
\textbf{Proof of Part 1: }
We compute $\psi( a_j a_i a_j ^ {-1} ) 	$.
\begin{align*}
	\psi (a_j a_i a_j ^{-1}) &= (1+x_j) (1+x_i) (1 - x_j + x_j ^2 - x_j ^3 \ldots ) \\
 &= 1 + ( x_i - x_i x_j + x_i x_j ^2 \ldots) + ( x_j x_i -x_j x_i x_j + x_i x_j ^2 \ldots)
 \end{align*}
 Note that the term $ x_i $ has coefficient $ 1 $. Since $ \lbrace i \rbrace $ is a proper
 subset of $ \lbrace j, i \rbrace $ and $ \lbrace j, i \rbrace $, it follows that
 $ \psi ( a_j a_i a_j ^{-1} ) \equiv \psi (a_i) \text{ mod } D_{i,n} $.
 It immediately follows that
\begin{equation*}
 \psi ( w_i ) \equiv \psi ( a_j w_i a_j ^{-1} ). \qed 
 \end{equation*}

\textbf{Proof of Part 2:}
We now demonstrate that $ \psi ( [a_j, w_j ] ) = 1 $. 
From part 1, we note that
\begin{align*}
\psi ( a_j w_j a_j ^{-1} ) &= \psi (w_j).  
\end{align*}
Now, 
\begin{align*}
\psi ( [a_j, w_j ] ) &= \psi ( w_j (w_j)^{-1} ) \\
&= \psi (1). \qed
\end{align*} 

\textbf{Proof of Part 3:}
See Lemma \ref{an}. \qed

\begin{thm} \label{fundth} The term $ \bar{ \mu} (x_{i_1}, x_{i_2}, \dots x_{i_s},w_i ) $ is a link invariant. \end{thm}
\textbf{Proof:} We note that the fundamental group of a virtual (welded) link diagram is invariant under the Reidemeister moves and the extended Reidemeister moves. To show that this is a link  invariant, we need only apply Lemma \ref{linkinvar}. This demonstrates that $ \bar{ \mu} $ is not altered by a change of base point or multiplication by elements equivalent to $ 1 $. \qed

We now consider the effect of crossing changes on the $ \bar{ \mu } $ invariants.
To demonstrate that $ \bar{ \mu } $ is a homotopy invariant, we must prove that the image of $ w_i $ in $ \mathbb{Z} [[x_1, x_2, \ldots x_n]] / D_{i,n} $ is invariant under crossing change. We consider the effect of a self crossing change on $ w_i $. 

Suppose that a crossing change occurs on the $ j^{th} $ component. Locally, the meridians of the $j^{th} $ arc are replaced by a different conjugate of $ a_j $. Let 
$ w_i = b_1 b_2 \ldots b_m $. Let $ b_1 = c a_j c $ (where $ c $ is a conjugate of $ a_{j} $, indicating a self-crossing). The term $ b_1 $ represents the $i^{th} $ component under passing the $ j^{th} $ strand. After a self-crossing change, $ w_i = \hat{b}_1 b_2 \ldots b_n $ and $ \hat{b}_1 = \hat{c} a_j \hat{c}^{-1} $ where $ \hat{c} $ is a conjugate of $ a_j $. To replace 
$ b_1 $ with $ \hat{b}_1 $, we multiply $ w_i $ by $ \hat{c} a_{j} \hat{c} ^{-1} a_{j}^{-1} 
a_j c a_j ^{-1} c^{-1} $. This term is an element of the \textit{commutator subgroup generated by $ a_j $} that consists of
elements of the form: $[k, \hat{k}]  $ where both $ k $ and $ \hat{k} $ are conjugates of $ a_j $ . We denote this subgroup as $[a_j]_2 $.

We now consider a crossing change on the $ i^{th} $ strand. Let $ w_i = b_1 b_2 \ldots b_m $ and suppose that $ b_1 $ represents the $i^{th} $ strand under passing itself, that is $ b_1 = c a_i^{ \pm 1} c^{-1} $. After the crossing change, $ w_i = b_2 \ldots \hat{b}_k \ldots b_m $, where $ \hat{b}_k  $ is the new over crossing (generated by the crossing change). Note that $ \hat{b}_k  = \hat{c} b_1^{-1} \hat{c}^{-1} $ , and to replace $ b_1 $ with $ \hat{b}_k $,we multiply 
$ w_i $ by $
( b_2 \ldots b_{k-1} \hat{c} b_1 ^{-1} \hat{c}^{-1} b_{k-1} ^{-1} \ldots b_2 ^{-1} ) b_1 ^{-1} $, an element of the subgroup generated by conjugates of $ a_i $.
Therefore, we need only show that $w_i $ is invariant under multiplication by an element of $ \langle a_i \rangle $, the subgroup generated by conjugates of $ a_i $.

\begin{lem}	\label{cut}
Let $ L $ be a $k $ component link, and let $w_i $ represent the longitude of the $i^{th} $ component. 
\begin{enumerate}		
	\item The Magnus expansion of $w_i$ is unchanged when multiplied by an element of the normal subgroup $[ a_{j}]_2 $ (a crossing change on the $ j^{th} $ strand).
	\item The Magnus expansion of $w_i$ is unchanged when multiplied by an element of the group $ \langle a_i \rangle $, the generator of the meridian of the $i^{th} $ component. (a crossing change on the $i^{th} $ strand).
\end{enumerate}
\end{lem}

\textbf{Proof of Lemma \ref{cut}, Part 1:}
Suppose $ X $ is an element of $ [a_j]_2$. Then we may assume that
$ X=[c a_j c^{-1} ,a_j^{-1} ] = [c, a_j] [c, a_j ^{-1} ] $.
Note that $[c, a_j] $ and $ [c, a_j ^{-1} ] $ are elements of $ A_2 $.  Hence, each non-constant term in $ \psi ([c,a_j]) $ and $ \psi ([c,a_j ^{-1}]) $ has degree 2 or greater and each term has contains at least two copies of $x_j$ since $ \psi(c)= 1+x_j \ldots $. Now, each non-constant term of $ \psi(X) $ has degree 2 or higher and contains at least two copies of $x_j $. However, since $ \mu (J,w_i) \equiv 0 $ if $ J$ contains duplicate entries, we see that $ \psi ( X w_i) \equiv \psi (w_i) $. \qed

\textbf{Proof of Lemma \ref{cut}, Part 2:}
Let $ \psi (w_i) = 1+ W $.
We may assume that $ y = c a_i c^{-1} $. Applying Lemma Let $ \psi (y) = 1 + Y $, where each term in $Y$ has at least one copy of $x_i $, so that $1+ Y \equiv 1 \in D_{i,n} $. 
Now, $ \psi(y w_i) = (1+Y)(1+W)$ so that
\begin{gather*}
\psi (y w_i) = 1 + Y + W + Y W \\
\equiv 1+ W \text{ modulo } D_{i,n} \qed
\end{gather*} 
We apply this lemma to prove that
 $ \bar{ \mu} (x_{i_1}, x_{i_2}, \ldots x_{i_s},w_i) $ is a link homotopy invariant.
 	
\begin{thm} The term $ \bar{ \mu} (x_{i_1}, x_{i_2}, \dots x_{i_s},w_i ) $ is a homotopy link invariant. \end{thm}
\textbf{Proof:} We apply lemma \ref{cut} to demonstrate that 

$\bar{\mu }(x_{i_1}, x_{i_2}, \ldots x_{i_s},w_i) $, obtained from the Magnus expansion of the reduced representative of the $i^{th} $ longitude, is unchanged by a self crossing change. Combined with theorem \ref{fundth}, this demonstrates that $ \bar{ \mu} $ is a link homotopy invariant. \qed

\section{Computational Examples} \label{examples}

In this section, we compute some simple examples using the Magnus expansion. We then apply an analog of  Polyak's skein relation \cite{polyak} for the $ \mu $ invariant to the same examples. Our construction of Polyak's skein relation for the $ \mu $ invariants utilizes the Magnus expansion and can be applied to virtual braids. Since the $ \bar{ \mu } $ invariants are computed progressively from the $ \mu $ invariants, this effectively allows us to compute the $ \bar{ \mu } $ invariants for the closure of a virtual braid. 
We give our own proof, based on the Wirtinger presentation, so that it is clear that this relation applies to virtual links and knots. The proof of this skein relation relies on the Magnus expansion and our earlier construction, so we begin by introducing some notation and proving several lemmas about infinite series. 
Let 
\begin{equation} \label{barnote}
\bar{v} = \sum_{i=1} ^{ \infty} (-1)^i v^i 
\end{equation}
then
\begin{equation*}
\bar{w} x w = \sum_{i=1} ^{ \infty} (-1)^i w^i xw.
\end{equation*}

The following lemmas will simplify computations in our proof. 
\begin{lem} \label{lem1} Let $ \bar{w} $ be defined as in equation \ref{barnote}. The equations: $ (1+w) (1 + \bar{w} ) =1$ and $ (1 + \bar{w})(1+ w) =1 $ both hold.
\end{lem}
\textbf{Proof:}
We observe that $ (1 + w)(1+ \bar{w}) $ expands to $ 1 + w + \bar{w} + w \bar{w} $. Expanding the individual terms, we observe that this expression is equivalent to one. \qed

\begin{lem} \label{lem2} The expression: $ (1 + \bar{w}) (1 + x) (1 + w) $ is equivalent to $ 1 + x + \bar{w} x + xw + \bar{w} x w $.
\end{lem}
\textbf{Proof:}
We expand $ (1 + \bar{w}) (1 + x) (1 + w) $ :
\begin{gather*}
= (1+x + \bar{w} + \bar{w} x ) (1+w) \\
= 1 + x + \bar{w} + \bar{w} x + w + xw  + \bar{w} w + \bar{w} x w .
\end{gather*}
This reduces to:
\begin{equation*}
1 + x + \bar{w} x + xw + \bar{w} x w. \qed
\end{equation*}

\begin{lem} \label{lem3} The expression: $ (1 + w) (1 + \bar{x}) (1 + \bar{w}) $ is equivalent to $ 1 + \bar{x} + w \bar{x}  + w \bar{x} \bar{ w} $.
\end{lem}
\textbf{Proof:}
The proof is similar to the proof of Lemma \ref{lem2}. \qed

\begin{figure}[htb] \epsfysize = 2.0 in
\centerline{\epsffile{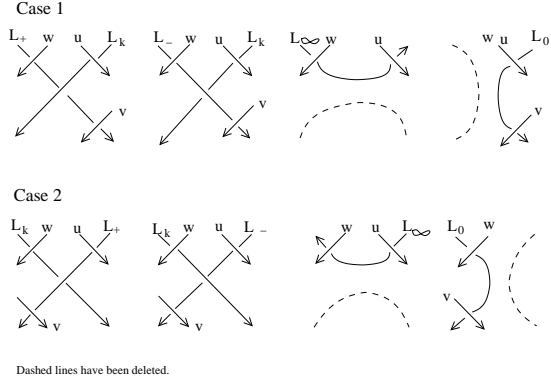}}
\caption{Polyak's skein relation}
\label{fig:pskein2}
\end{figure}

\begin{thm} The $ \mu $ invariant of the braids shown in figure \ref{fig:pskein2} can be computed based on the crossing between the pairs: $L_+ $, $L_k$ and $L_-$, $L_k$. Note that figure \ref{fig:pskein2} contains information that allows us to compute the longitudes. Let $ l_+ , l_{-}, l_{0},$ and $ l_{ \infty} $ denote the longitude of the the strands labeled $L_+, L_{-}, L_0, $ and $L_{ \infty} $. Let $l_k$ denote the initial label on the $k^{th} $ strand ($L_k$) so that $ \psi (l_k) = 1 + x_k $. Then the following skein relations hold for the $ \mu $ invariants. In case 1,
\begin{align*}
\mu (x_k, l_{+}) - \mu(x_k, l_{-}) &= 1  \\
\mu (x_{i_1}, x_{i_2}, \ldots x_{i_n},x_k, l_{+}) - \mu(x_{i_1}, x_{i_2}, \ldots x_{i_n},x_k, l_{-}) &= \mu ( x_{i_1}, x_{i_2}, \ldots x_{i_n}, l_0)  \\
\mu (x_k, x_{i_1}, x_{i_2}, \ldots x_{i_n}, l_{+}) - \mu(x_k, x_{i_1}, x_{i_2}, \ldots x_{i_n}, l_{-}) &= \mu ( x_k, x_{i_1}, x_{i_2}, \ldots x_{i_n}, l_{\infty} )  \end{align*}
and
\begin{gather*}
\mu ( x_{i_1}, x_{i_2}, \ldots x_{i_j}, x_k, x_{i_{j+1}} \ldots x_{i_n}, l_{+}) - \mu(x_k, x_{i_1}, x_{i_2}, \ldots x_{i_j} x_k, x_{i_{j+1}} \ldots x_{i_n}, l_{-}) \\ = \mu ( x_{i_1}, x_{i_2}, \ldots x_{i_j}, l_{0} )  \mu (x_{i_{j+1}} \ldots x_{i_n}, l_{\infty}).
\end{gather*}
In case 2, 
\begin{align*}
\mu (x_k, l_{+}) - \mu(x_k, l_{-}) &=1  \\
\mu (x_{i_1}, x_{i_2}, \ldots x_{i_n},x_k, l_{+}) - \mu(x_{i_1}, x_{i_2}, \ldots x_{i_n},x_k, l_{-}) &= \mu ( x_{i_-}, x_{i_2}, \ldots x_{i_n}, l_0)  \\
\mu (x_k, x_{i_1}, x_{i_2}, \ldots x_{i_n}, l_{+}) - \mu(x_k, x_{i_1}, x_{i_2}, \ldots x_{i_n}, l_{-}) &= \mu ( x_k, x_{i_1}, x_{i_2}, \ldots x_{i_n}, l_{\infty} )  
\end{align*}
\begin{gather*}
\mu ( x_{i_1}, x_{i_2}, \ldots x_{i_j}, x_k, x_{i_{j+1}} \ldots x_{i_n}, l_{+}) - \mu(x_k, x_{i_1}, x_{i_2}, \ldots x_{i_j} x_k, x_{i_{j+1}} \ldots x_{i_n}, l_{-}) \\ = \mu ( x_{i_1}, x_{i_2}, \ldots x_{i_j}, l_{0} )  \mu (x_{i_{j+1}} \ldots x_{i_n}, l_{\infty}).
\end{gather*}
\end{thm}
\textbf{Proof:} We compute the longitude of the components $ L_+, L_{-}, L_0,$ and $ L_{ \infty } $from figure \ref{fig:pskein2}. Based on the longitudes, we determine the Magnus expansion. Recall that the invariant  $ \mu (i_1, i_2, \ldots i_n , l) $ is the coefficient of the term $ x_{i_1} x_{i_2} \ldots x_{i_n} $ in $ \psi (l) $ (the Magnus expansion of the longitude $l$). We consider only case 1, case 2 is analogous.
The longitude of $L_+ $ is:
\begin{equation} \label{l+}
l_+ = vu^{-1} l_k u w.
\end{equation}
The longitude of $L_-$ is:
\begin{equation} \label{l-}
l_- = v w.
\end{equation}
The longitude of $L_0$ is:
\begin{equation} \label{l0}
l_0 = v u^{-1}.
\end{equation}
The longitude of $L_{\infty} $ is:
\begin{equation} \label{linf}
l_{ \infty} = u w.
\end{equation}
Note that $w$, $u$, and $v$ are words in the Chen group. The term $l_k$ is the initial label on the $k^{th} $ strand, so that the label on this strand when it crosses $L_+ $ is $ w ^{-1} l_k w $. Let $ \psi(v) = 1 + V$, $ \psi (w) = 1 + W$, 
$ \psi (u) = 1 + U $, and $ \psi (l_k) = 1 + x_k $. We observe that $ \psi( w^{-1} ) = 1 + \bar{W} $.
Now,
\begin{equation} \label{poseq}
\psi (l_+) = (1 + V)(1+ \bar{U} ) (1 + x_k) (1 + U) (1 + W) \end{equation}
Expanding, we determine that
\begin{gather*}
\psi (l_+)
= 1 + V + x_k + V x_K + \bar{U} x_k + V \bar{U} x_k \\ +  x_k U + V x_K U + \bar{U} x_k U + V \bar{U} x_k U \\
W + VW + x_k W + V x_k W+ \bar{U} x_k W + V \bar{U} x_k W + x_k U W\\ + V x_K U W + \bar{U} x_k U W + V \bar{U} x_k U W. 
\end{gather*}
We compute the value of $ \psi (l_-) $ based on Lemmas \ref{lem1} and \ref{lem2}
\begin{equation} \label{negeq}
\psi(l_-) = (1+V)(1+W) = 1 + V + W + VW.
\end{equation}
Next, we compute the $L_0 $ and $L_{ \infty} $ terms
\begin{equation} \label{infeq}
\psi( l_{ \infty} ) = (1 + U)(1 + W) = 1 + U + W + UW.
\end{equation}
and
\begin{equation} \label{0eq}
\psi(l_0) = (1 + V)(1 + \bar{U}) = 1 + V + \bar{U} + V \bar{U} .
\end{equation}

We consider $ \mu (x_k, l_+) $. This is the coefficient of the term $ x_k $ in $ \psi (l_+ ) $ (equation \ref{poseq}). We observe that $ \psi (l_-) $ (equation \ref{negeq}) does not contain an $x_k $ term, hence $ \mu (x_k, l_-) = 0$. We observe that the coefficient of $x_k $ in 
$ \psi (l_+) $ (equation \ref{poseq}) is one. As a result,
\begin{equation*}
\mu (x_k, l_{+}) - \mu(x_k, l_{-}) = 1.
\end{equation*}
We next consider $ \mu (x_{i_1}, x_{i_2}, \ldots x_{i_n}, x_k, l_+) $ and  $ \mu (x_{i_1}, x_{i_2}, \ldots x_{i_n}, x_k, l_-) $ which are the coefficients of monomials of the form  $ x_{i_1}, x_{i_2}, \ldots x_{i_n}, x_k $ in $ \psi(l_+) $ and 
$\psi (l_-) $ respectively. 
We note that $ \psi(l_- ) $ contains no monomials that end in $ x_k$ and $ \mu (x_{i_1}, x_{i_2}, \ldots x_{i_n}, x_k, l_+) =0 $. Monomials (that do not have the form $ x_k$) in $ \psi (l_+) $ with the last term $ x_k $ are: $ \bar{U}x_k, vx_k,$ and
$ V \bar{u} x_k $. These monomials correspond to the monomials in $ \psi (l_0) $ (equation \ref{0eq}).
Therefore:
\begin{equation*}
\mu (x_{i_1}, x_{i_2}, \ldots x_{i_n},x_k, l_{+}) - \mu(x_{i_1}, x_{i_2}, \ldots x_{i_n},x_k, l_{-})  = \mu ( x_{i_1}, x_{i_2}, \ldots x_{i_n}, l_0). 
\end{equation*}
We find all monomials in $ \psi (l_+) $ (equation \ref{poseq}) not equal to $x_k$ where the first term is $x_k$:
$x_k W, x_k U,$ and $x_k UW$. These terms correspond to the monomials in $ \psi (l_{\infty}) $. Hence:
\begin{equation*}
\mu (x_k, x_{i_1}, x_{i_2}, \ldots x_{i_n}, l_{+}) - \mu(x_k, x_{i_1}, x_{i_2}, \ldots x_{i_n}, l_{-}) = \mu ( x_k, x_{i_1}, x_{i_2}, \ldots x_{i_n}, l_{\infty} ).
\end{equation*}
We consider monomials (not equal to $x_k $) where $x_k$ occurs as an intermediate term. In $ \psi(l_-) $, we have the terms:
\begin{gather*}
V x_k U, \text{ } \bar{U} x_k U,\text{ }  V \bar{U} x_k U,\text{ }  V x_k U W, \text{ } \bar{U} x_k U W, \text{ and }V \bar{U} x_k U W. 
\end{gather*}
These correspond to the monmials in the product of $ \psi (l_0) $ and $ \psi (l_\infty) $.
We conclude that:
\begin{gather*}
\mu ( x_{i_1}, x_{i_2}, \ldots x_{i_j}, x_k, x_{i_{j+1}} \ldots x_{i_n}, l_{+}) - \mu(x_k, x_{i_1}, x_{i_2}, \ldots x_{i_j} x_k, x_{i_{j+1}} \ldots x_{i_n}, l_{-}) \\ = \mu ( x_{i_1}, x_{i_2}, \ldots x_{i_j}, l_{0} )  \mu (x_{i_{j+1}} \ldots x_{i_n}, l_{\infty}). \qed
\end{gather*}

\begin{rem} Although we can use Polyak's skein relation to compute the $ \mu $ invariant, we can not use his result to directly determine the $ \bar{ \mu } $ invariants. To compute the $ \bar{ \mu } $ invariants we need to reduce modulo $ D_{i,n} $. We should also note that the closure of two inequivalent braids may result in two equivalent welded (virtual) braids.
\end{rem}

\subsection{Virtual Hopf Link}

\begin{figure}[htb] \epsfysize = 2.0 in
\centerline{\epsffile{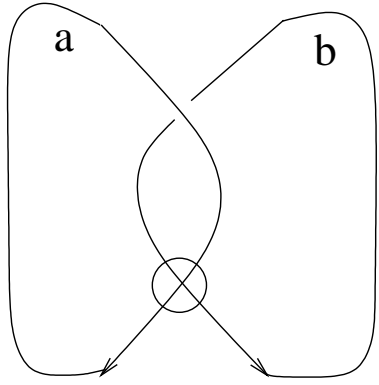}}
\caption{Braid whose closure is the virtual Hopf link}
\label{fig:vhopf}
\end{figure}

The closure of the braid shown in figure \ref{fig:vhopf} is the virtual Hopf link. 

We compute the fundamental group of the virtual Hopf link is generated by $a $ and $ b $ modulo the following relation:
\begin{equation*}
 b = a b a^{-1}
 \end{equation*}
The longitude of the strands is given by:
 \begin{align*} 
 w_a &= 1 \\
 w_b & = a^{-1}
 \end{align*}

We apply the Magnus expansion (sending $ b \rightarrow 1+b $ and $ a \rightarrow 1+a $ ) to the longitudes:
 \begin{equation*}
 \psi (a) = 1- a +a^2 \ldots \text{ and } \psi (1) =1. 
 \end{equation*}
There are no proper subsets of $ \lbrace b \rbrace $ and $ \lbrace a \rbrace $. 
As a result, $ \bar{ \mu } (b,w_a) = 0 $ and $ \bar{ \mu } (a, w_b) =-1 $. 
Note that $ link(b,a) =0 $ and that $ link(a,b) = -1 $.

We apply case 2 of the skein relation to the classical crossing the virtual Hopf link and observe that it is correct.
We note that $ \mu (a, w_b) = -1 $ which corresponds to $ \mu (x_k, l_-)$. Switching the crossing, we obtain
$ \mu (a, w_{b+}) =0 $  which corresponds to $ \mu (x_k, l_+) $.
\begin{equation*}
\mu ( a,  w_{b+}  ) - \mu (a, w_{b}) = 1.
\end{equation*}

\subsection{Virtual Link Example}

\begin{figure}[htb] \epsfysize = 2.0 in
\centerline{\epsffile{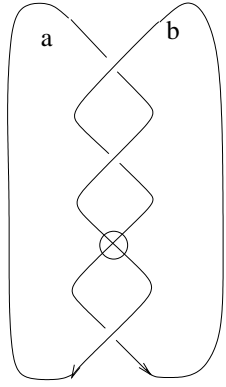}}
\caption{Virtual link example}
\label{fig:virtexa}
\end{figure}

The closure of the of the braid shown in figure \ref{fig:virtexa} is a two component virtual link.

Observe that $ link(b,a) = 1 $ and $ link(a,b) = 2 $.
The fundamental group of the link is generated by $ a $ and $ b $. 
The relations in fundamental group are:

\begin{align*}
a &= b a b ^{-1} \\
b &=  a ^2 b a ^ {-2}.
\end{align*}
The longitudes of the strands are:
 $ w_a = b $ and $ w_b = a ^2 $.

Applying the Magnus expansion to the reduced longitudes:
\begin{align*}
\psi (w_a) &= 1+b \\
\psi (w_b) &=  (1+a)^2 \\
&= 1 + 2 a +  a^2.
\end{align*}
Since there are two strands, all elements with degree greater than 1 are in $ D_{a,2} $ and $ D_{b,2} $.
Now, $ \mu(b,w_a) = \bar{ \mu} (b,w_a) =1 $ and $ \mu (a, w_b) = \bar{ \mu }(a, w_b)  = 2 $. 

We apply the skein relation to the uppermost crossing in this braid. 
We observe that $ \mu (a, w_b) = 2 $ and that for the switched crossing: $ \mu (a , w_{b-} ) =1 $
$ \mu (a,w_b) - \mu(a,w_{b-}) = 2-1 $. Next, $ \mu (b,w_a) = 0 $ and switching the crossing, we obtain $ \mu (b, w_{a+} ) = 1$. Hence $ \mu (b, w_{a+}) - \mu (b,w_a) = 1-0$.

\subsection{Modification of the Borromean Links}

\begin{figure}[htb] \epsfysize = 2.0 in
\centerline{\epsffile{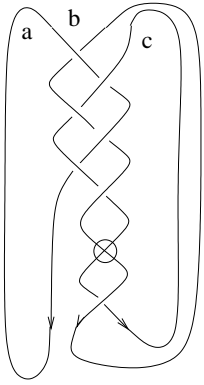}}
\caption{Modification of the Borromean links}
\label{fig:modboro}
\end{figure}

The closure of the braid shown in figure \ref{fig:modboro} is obtained by adding a classical and virtual twist to the Borromean links on  $ 2^{nd} $ and $ 3^{rd} $ strands.
If the $2^{nd} $ or $ 3^{rd} $ strand is removed from the closure of this braid, we obtain an unknotted 2 component unlink. If the first strand is removed, we obtain a 2 component link.

The fundamental group of this link is generated by $ a, b, c $. We obtain the following relations:
\begin{gather*} 
a =a ^{-1} b ^{-1} a c ^{-1} a ^{-1} b a c a c ^{-1} a ^{-1} b ^{-1} a c a ^{-1} b a \\
b = c a c ^{1} a^{-1} b a c a ^{-1} c ^{-1} \\
c = c a c ^{-1} a^{-1} b ^2 a c a ^{-1} c ^{-1} a ^{-1} b^{-1} a  c   a^{-1} b a c a c^{-1} a^{-1} b ^{-2} a c a ^{-1} c ^{-1} .
\end{gather*}
The longitude of these components:
\begin{gather*}
w_a = a ^{-1} b ^{-1} a c ^{-1}  a ^{-1} b a c \\
w_b = c a c ^{-1} a ^{-1} \\
w_c =  c a c ^{-1} a ^{-1} b ^{2}  a c a ^{-1} c ^{-1} a ^{-1} b ^{-1} a .
\end{gather*}
We reduce the longitudes modulo the second communtator group: 
\begin{gather*} 
w_a \equiv ( a ^{-1} b ^{-1} a  c ^{-1}  a ^{-1} b a)( c ) \\
\equiv  b ^{-1} c^{-1}  b  c \text{ modulo } A_2 \\
w_b \equiv c a c^{-1} a^{-1} \\
w_c \equiv ( c a c ^{-1} a ^{-1} b ^{2}  a c a ^{-1} c ^{-1}) ( a ^{-1} b ^{-1} a ) \\
\equiv b^2 a^{-1} b^{-1} a 
\end{gather*}
We now apply the Magnus expansion to the reduced form of $w_a$, that is: $ b^{-1} c^{-1} b c $.
\begin{gather*}
\psi ( b^{-1}  c^{-1}  b  c ) = (1-b + b^2 \ldots ) ( 1-c + c^2 \ldots ) (1+b) (1+c) \\
= (1-b + b^2 \ldots -c + b c - b^2 c \ldots + c^2 - b c^2 + b^2 c^2 \ldots) 
( 1+b + c +bc) \\
= 1-b + b^2 \ldots -c + b c - b^2 c \ldots 
+ c^2 - b c^2 + b^2 bc^2 \ldots \\ 
+ b - b^2 -cb 
+ c - b c \ldots + b c ) 
\end{gather*}
All terms with degree greater than 2 are in $ D_{a,3} $.
Collecting coefficients, we observe that:
\begin{align*}
\mu ( b, w_a) &= 0 \\
\mu (c , w_a) &= 0 \\
\mu (b,c, w_a) &= 1 \\
\mu (c,b, w_a) &= -1
\end{align*}

The computation for $ \psi (w_b) $ is similar. Hence:
\begin{align*}
\mu ( a, w_b) &= 0 \\
\mu (c , w_b) &= 0 \\
\mu (a,c, w_b) &= 1 \\
\mu (c,a, w_b) &= -1
\end{align*}

We compute $ \psi (w_c) = \psi (b^2 a^{-1} b^{-1} a ) $.
Note that
\begin{gather*}
(1+b)^2  (1-a + a^2 \ldots) (1-b+b^2 \ldots) (1+a)\\
 = 1+b + ab -ba - b^2 \ldots
\end{gather*}
From this equation, we obtain:

\begin{align*}
\mu (a,w_c) &= 0 \\
\mu (b,w_c) &= 1 \\
\mu (a,b,w_c) &= 1 \\
\mu (b,a, w_c) &= -1
\end{align*}

Note that since $ \mu (i, w_j) =0 $ for all $ i,j $ except for the case where $i=b $ and $ j=c $, that $ \mu (i,k, w_j) = \bar{ \mu } (i,k,w_j) $.

\subsection{A Virtual Link with Linking Number Zero} \label{lastexa}
\begin{figure}[htb] \epsfysize = 2.0 in
\centerline{\epsffile{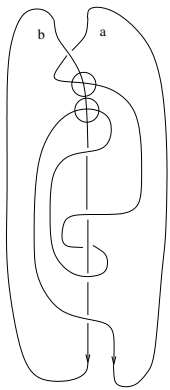}}
\caption{Modification of a link homotopic to the identity}
\label{fig:homotopyexa}
\end{figure}

The closure of the link shown in figure \ref{fig:homotopyexa} is obtained by changing two classical crossings in a link homotopic to the identity into virtual crossings. 
\begin{rem} Note that $ link(a,b)= 0 $ and $link(b,a) =0 $ for this link. A classical, 2 component link with linking number zero is homotopically trivial. If this link is homotopic to a classical link then it is homotopic to the unknotted, 2 component unlink.
\end{rem} 

The fundamental group is generated by $ a $ and $ b $ with the relations:
\begin{align*}
a &=b^{-1} a b a^{-1} b a b a^{-1} b^{-2} a b^2 a b^{-1} a^{-1} b^{-1} a b a^{-1} b \\
b &= a b a^{-1} b a b a^{-1}.
\end{align*}
From these relations, the longitudes of the components are:
\begin{align*}
w_a &= b^{-1} a b a^{-1} b a b a^{-1} b^{-2} \\
w_b &= a b a^{-1}.
\end{align*}
We apply the Magnus expansion to $ w_a $:
\begin{equation*}
\psi (w_a) \equiv 1 \text{ mod } D_{a,2}
\end{equation*}
Applying the Magnus expansion to $w_b$:
\begin{equation*}
\psi (w_b) = 1 + (b - ba + ba^2 - ba^3 \dots) + (ab -aba + a ba^2 - a b a^3 \dots) 
\end{equation*}
We conclude that 
\begin{equation*}
\mu ( b, w_b) = 1.
\end{equation*}
Note that $ \bar{ \mu } (b, w_b) \equiv 0 $ by definition.
In the next section, we demonstrate that $ \mu (i,w_i ) =0 $ for all components of a classical link and a link homotopic to a classical link, allowing us to conclude that this link is not homotopic to any classical link. 

\section{The $ \mu $ invariant and virtual link diagrams} \label{links}

\begin{figure}[htb] \epsfysize = 0.5 in
\centerline{\epsffile{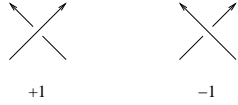}}
\caption{Crossing sign}
\label{fig:sign}
\end{figure}

In this section, we consider the invariant $ \mu (a, w_b) $ for $ a \neq b $ and $ a=b $. 
We recall that the definition of \emph{linking number} for virtual link diagrams. 
The sign of a crossing, $ c $ is determined by the local orientation as shown in figure \ref{fig:sign}.

We define the \emph{link(b,a)} to be:
\begin{equation}
\sum_c sign(c)
\end{equation}
where $ c $ is a crossing such that the strand $ a $ passes under strand $b $. 
Note that $ link(b,a) $ is not necessarily equivalent to $link (a,b)$.
\begin{rem} In the classical case, $ link(a,b) = link(b,a) $ for $ a=b $ and $ a \neq b $. 
\end{rem}
\begin{rem} In \cite{milnor2}, Milnor observes that $ \mu (j,k,w_i) = \mu (k,j,w_i) $ and develops the idea of shuffles, symmetry relations on the $ \bar{ \mu} $ invariants. The symmetry relations need to be re-examined for virtual and welded diagrams.
\end{rem}
We will demonstrate that $ \bar{\mu} (b,w_a) = \mu (b,w_a) = link (b,a) $ for $ a \neq 0 $.  
\begin{thm}Let $ L $ be a link with at least two components, $ a $ and $ b $. Then $ \bar{ \mu} (b,w_a) = \mu (b, w_a) = link(b,a) $.
\end{thm} 
\textbf{Proof:} Consider the strand $a$. If strand passes under strand $ b $ with a positive orientation $k $ times and with a negative orientation $j$ times, then $ link(b,a) = k-j $.

We compute the longitude of $ a $, $ w_a $. When $a_i $ undercrosses $ b_j $, we note that $ a_{i+1} = (b^c) a_i (b^c) ^{ -1} $, where $ b^c $ is a conjugate of $ b $. Hence, for each positive undercrossing, $ w_a $ contains a $ b $ and for each negative crossing a $ b^{-1} $. 
Now, applying the Magnus expansion to $ w_a $, we see that $ k $ factors have the form $ (1+b)$ and $ j $ factors have the form $ (1-b + b^2 -b^3 \ldots ) $. All other terms are of the form $ (1+c) $ or $ (1 -c +c^2 -c^3 \ldots ) $. To compute the coefficient of $ b $ in the Magnus expansion,
\begin{align}
(1+b)^k (1-b+b^2-b^3 \ldots )^j &= (1+ k b + \ldots) (1- j b + \ldots) \\
&= 1+ (k-j) b + \ldots
\end{align}
Hence $ \mu (b,w_a) = k-j $. We observe that $ \lbrace b \rbrace $ does not contain a proper subset. Then $ \bar{ \mu} (b,w_a) = k-j $.  \qed

We now consider $ \mu (a, w_a) $. In the classical case, $ \mu (a,w_a) =0 $ for all components. For the example given in section \ref{lastexa}, $ \mu (b, w_b) =1 $.
\begin{thm} Let $ L $ be a classical $n$ component link diagram. Then $ \mu (i,w_i ) $ $= 0 $ for any component $i$ of the braid whose closure is $L$.  
\end{thm}
\textbf{Proof: } Let $ L $ be a classical $ n $ component link diagram. Then $ L $ is homotopic to the closure of a pure classical braid and does not have any self-crossings. Consider the generator of the meridian of component $ i $, $a_i $. The longitude of this component, $w_i$ has the form: $ w_i =l_1 l_2 \ldots l_n  $ where each $ l_k $ is conjugate of the generator of a component that underpasses component $ i $. Hence,
 $ l_k $ is not a conjugate of $ a_i $ for all $ k \in \lbrace 1,2, \dots n \rbrace $. The meridian $ a_i $ occurs only as a part of a conjugate pair in the longitude. Hence, $ \mu (i,w_i) = 0 $.\qed

\begin{cor} Let $ L $ be a link diagram. Let $ i $ be a component of $L $ with longitude $w_i $. If $ \mu (i, w_i) \neq 0 $ then $ L $ is not a classical link diagram. 
\end{cor}
\textbf{Proof: }The proof of this corollary is immediate.\qed
\begin{rem}
We can amend Polyak's skein relation to include $ \mu (i,w_i) $. Let $ L $ be an $ n $ component diagram and let $ w_i ^+ $ represent the $i^{th} $ component with a fixed positive self crossing. Let $ w_i ^- $ represent the diagram with the fixed self crossing changed to a negative crossing. Then
\begin{equation*}
\mu (i,w_i ^+) - \mu (i, w_i ^-) \equiv 0 \text{ modulo } 2.
\end{equation*}
\end{rem}

\section*{Acknowledgments}The authors would like to thank Xiao Song Lin for the suggestion to investigate virtual link homotopy.
This research was performed while the first author held a National Research Council Research Associateship Award jointly at the Army Research Laboratory and the United States Military Academy. 
The U.S. Government is authorized to reproduce and distribute
reprints
for Government purposes notwithstanding any copyright annotations thereon.
The
views and conclusions contained herein are those of the authors and should
not be
interpreted as necessarily representing the official policies or
endorsements,
either expressed or implied, of the United States Military Academy or the U.S. Government. It gives the second author great pleasure to acknowledge support from  NSF Grant DMS-0245588.

\end{document}